\documentclass[12pt]{amsart}
\usepackage{amsmath,amssymb,amsfonts,amscd,graphicx}

\usepackage{framed}

\usepackage{soul}

\usepackage{xcolor}

\usepackage{multirow}

\newtheorem{theorem}[subsection]{Theorem}

\newtheorem{proposition}[subsection]{Proposition}

\theoremstyle{definition}
\newtheorem{definition}[subsection]{Definition}
\newtheorem{conjecture}[subsection]{Conjecture}

\theoremstyle{remark}

\numberwithin{equation}{section}

\theoremstyle{definition}

\newcommand{\co}{\colon\thinspace}

\begin{document}

\title[Geometric complexity of embeddings in ${\mathbb R}^d$]{Geometric complexity of embeddings in ${\mathbf{\mathbb R}^d}$}
\author[Michael Freedman and Slava Krushkal]{Michael Freedman and Vyacheslav Krushkal}
\address{Microsoft Station Q, University of California, Santa Barbara, CA 93106-6105}
\email{michaelf\char 64 microsoft.com}
\address{Department of Mathematics, University of Virginia, Charlottesville, VA 22904}
\email{krushkal\char 64 virginia.edu}

\thanks{VK is partially supported by the NSF}

\dedicatory{Dedicated to Misha Gromov on the occasion of his 70th birthday}

\begin{abstract}
Given a simplicial complex $K$, we consider several notions of geometric complexity of embeddings of $K$ in a Euclidean space ${\mathbb R}^d$:
thickness, distortion, and refinement complexity (the minimal number of simplices needed for a PL embedding).
We show that any $n$-complex with $N$ simplices which topologically embeds in ${\mathbb R}^{2n}$, $n>2$,
can be PL embedded in ${\mathbb R}^{2n}$ with refinement complexity
$O(e^{N^{4+{\epsilon}}})$.
Families of simplicial $n$-complexes $K$ are constructed such that any embedding of $K$ into ${\mathbb R}^{2n}$ has 
an exponential lower bound on thickness and refinement complexity as a function of the number of simplices of $K$. 
This contrasts embeddings in the stable range, $K\subset {\mathbb R}^{2n+k}$, $k>0$, where all known bounds on geometric complexity 
functions are polynomial.
In addition, we give a geometric argument for a bound on distortion of expander graphs in Euclidean spaces. 
Several related open problems are discussed, including questions about the growth rate of complexity functions of embeddings, 
and about the crossing number and the ropelength of classical links.
\end{abstract}

\maketitle

\section{Introduction.}

We use the abbreviation $n$-complex to mean an $n$-dimensional simplicial complex.
By general position any $n$-complex embeds into the Euclidean space ${\mathbb R}^d$ in the stable range, $d\geq  2n+1$.
There have been recent developments, cf.  \cite{GG, Pardon}, in quantifying the geometry of embeddings, 
and polynomial ``complexity'' has been established in the stable range.  There are several possible notions of thickness 
and distortion of which we will mention four. 
We also consider refinement complexity, the minimal number of simplices needed for a PL 
embedding of a given complex into ${\mathbb R}^d$.
We prove that for the refinement complexity and for the notion of thickness defined below 
the lower bound shifts from polynomial 
to at least exponential when the dimension $d$ passes from $2n+1$ to $2n$. 
We do not presently know if there is a similar transition for the other notions
of geometric complexity of embeddings. 
For the refinement complexity we also establish an upper bound:
we prove that any $n$-complex, $n>2$, with $N$ simplices and with trivial van Kampen's obstruction
can be PL embedded in ${\mathbb R}^{2n}$  with refinement complexity
$O(e^{N^{4+{\epsilon}}})$. Our conclusions on thickness and refinement complexity 
are summarized in a table appearing near the end of section~\ref{n-complexes}.

Section \ref{relative section} 
considers the smallest codimensions $0, 1, 2$.
There are families of $(d-1)$- and $d$-complexes for which the refinement complexity  
with respect to embeddings into ${\mathbb R}^d$, $d\geq 5$,
is a non-recursive function of the number of simplices. We conjecture that the thickness (defined below) of these complexes is also non-recursive, see section \ref{relative section}.
This construction \cite{MTW}, stated as theorem \ref{non-recursive thm}, is based on 
Novikov's theorem on the algorithmic unrecognizability of spheres $S^n$,  $n\geq 5$, and
it  is motivated by the work \cite{MTW} on algorithmic undecidability of the embedding problem in codimensions $0,1$.

In theorem \ref{relative theorem} we show that in the relative setting there are also embedding problems in codimension $2$: pairs $(K,L)$ where $L$ is a subcomplex of a $(d-2)$-complex $K$ with a fixed embedding $L\subset {\mathbb R}^d$, $d\geq 5$,
such that the refinement complexity (and conjecturally thickness)  of any extension to an embedding $K\subset {\mathbb R}^d$
is a non-recursive function. These examples are based on groups with undecidable word problem.

The rest of the paper addresses
the question of what kinds of complexity functions can be realized in higher codimensions nearer the stable range. 
To formulate an upper bound for $n$-complexes in ${\mathbb R}^{2n}$ and to fix the notation we state the following definition.

\begin{definition} \label{rc}
{\bf Refinement complexity (rc)}. Given a complex $K^n$ and $d\geq n$, ${\rm rc}(K,d)$ is the minimal number of simplices 
in a PL embedding of $K$ into ${\mathbb R}^d$. We say ${\rm rc}(K,d)=\infty$ if there is no PL embedding $i\! : K^n\hookrightarrow {\mathbb R}^d$.
There is also an extrinsic form: given a piecewise mooth embedding $i\! : K\hookrightarrow {\mathbb R}^d$, ${\rm rc}(i)$ is the minimum number of simplices of a PL embeddings of $K$ piecewise smoothly isotopic to $i$.
\end{definition}

\begin{theorem} \label{rc upper bound thm} \sl
Let $K$ be an $n$-complex with $N$ simplices which topologically embeds in ${\mathbb R}^{2n}$, $n>2$.
Then ${\rm rc}(K, 2n)=O(e^{N^{4+{\epsilon}}})$ for any ${\epsilon}>0$.
\end{theorem}

The proof of this theorem in section \ref{rc bound}  is a quantitative version of the proof
that the vanishing of the van Kampen cohomological obstruction for $n$-complexes in ${\mathbb R}^{2n}$ is sufficient
for embeddability for $n>2$.

The main result of section \ref{n-complexes}
is the construction of $n$-complexes in ${\mathbb R}^{2n}$ with an (at least) exponential lower bound on refinement complexity and thickness. To formulate this result, we
first state two relevant notions of embedding thickness.

Consider embeddings  $K^n \longrightarrow B^d_1\subset {\mathbb R}^d$ of an $n$-dimensional simplicial complex into the unit ball in ${\mathbb R}^d$. 
Following \cite{GG}, we say
that such an embedding has {\em Gromov-Guth thickness} at least $T$ if the distance between the images of any two non-adjacent 
simplices is at least $T$. Assuming that local combinatorial complexity of $K$ is bounded, it is shown in \cite{GG} that 
if $d\geq 2n+1$ then  $N^{-\frac{1}{d-n}}$  is a sharp (up to an ${\epsilon}$ summand in the exponent) bound for the Gromov-Guth thickness of $K$ in ${\mathbb R}^d$, where $N$ is the number of vertices of $K$.
This is a generalization of earlier work of Kolmogorov and Barzdin \cite{KB} which, predating a formal definition of expander graphs, used their properties to establish the sharp bound $N^{-1/2}$ for embedding thickness of graphs in $3$-space.
Next we define the notion of thickness used in this paper.

\begin{definition} \label{defi thickness}
Consider embeddings $K^n \longrightarrow B^d_1\subset {\mathbb R}^d$ which are cell-wise smooth. 
(Explicitly this means that each closed $k$-simplex is carried into 
${\mathbb R}^d$ by some $C^{\infty}$ embedding of a neighborhood
${\mathcal N}({\sigma}^k)\subset {\mathbb R}^k$.)
We say that the {\em thickness} of such an embedding is at least $T$ if 
the distance between the images of any two non-adjacent 
simplices is at least $T$ and all simplices have embedded $T$-normal bundles. (In the codimension zero case
there is no normal bundle requirement on $d$-simplices embedded in ${\mathbb R}^d$.)
\end{definition}

This definition, like other notions of thickness and distortion, has two versions: thickness of an isotopy class of an embedding $i$ (where one takes the 
supremum over all embeddings isotopic to $i$), and {\em intrinsic thickness} of $K$ (the supremum over all embeddings into ${\mathbb R}^d$).

The notion of thickness introduced above, taking into account the size of the normal bundle, arises naturally when one considers isotopy classes of embeddings of manifolds.
In particular, related measures of complexity have been investigated in knot theory, cf. \cite{CKS, LSDR, Nabutovsky}.

In contrast to the known polynomial bounds on thickness, refinement complexity and distortion (discussed below) of embeddings in the stable range, we show the following.
\begin{theorem} \label{exponential theorem} \sl
For each $n\geq 2$ there exist families of $n$-complexes $\{ K_l\} $ with
$m_l\longrightarrow~\infty$ simplices and bounded local combinatorial complexity which embed into ${\mathbb R}^{2n}$, and the thickness of any such 
embedding is at most $c^{-m_l}$. 
Moreover,
$$ C^{m_l}\; <\; {\rm rc}(K_l, 2n)\; < \; \infty.$$
Here the constants $c, C>1$ depend on $n$.
\end{theorem}

{\bf Remarks}. 1. The phrase ``bounded local combinatorial complexity'' means that 
in a family $K_l$ the maximum number of simplices incident to any fixed simplex is bounded independently of $l$.

The proof of theorem \ref{exponential theorem} is given in section \ref{n-complexes}.
At the end of that section we state an addendum extending this theorem to $(n+k)$-complexes in ${\mathbb R}^{2n+k}$ for all $k\geq 0$.
A table summarizing our results on refinement complexity and thickness of $n$-complexes in ${\mathbb R}^d$ for various values of 
$(n,d)$ is included at the end of section \ref{n-complexes}.

Other related families of $2$-complexes in ${\mathbb R}^4$ are discussed in section \ref{2-complexes}. This construction
is a geometric implementation of nested commutators and other recursively defined elements in the free group. The key feature is that a word of length $2^k$ can be created using $O(k)$ $2$-cells. It seems likely that these examples also have exponentially small thickness as in theorem \ref{exponential theorem}, 
however at present this is an open problem. Section \ref{2-complexes} also formulates related questions about the thickness and the crossing number of links in ${\mathbb R}^3$.

Next we review some of the known results about distortion of embeddings into Euclidean spaces.
Given a subset $K\subset {\mathbb R}^n$ and two points $x,y\in K$, let $d_K (x,y)$ denote the intrinsic metric: the infimum of the lengths of paths in $K$ connecting
$x$ and $y$. The {\em distortion} of $K$ is defined as
\begin{equation} \label{distortion}
{\delta}(K)\; :=\; \sup_{x,y\in K}\; \frac{d_K(x,y)}{d_{{\mathbb R}^d}(x,y)}.
\end{equation}

The distortion of an isotopy class of an embedding $i\co K\longrightarrow {\mathbb R}^d$  is the 
infimum over all embeddings isotopic to $i$. Measuring the distortion is often a rather subtle problem.
A recent advance in the subject is due to J. Pardon \cite{Pardon} who used integral geometry of curves on surfaces to show 
that the distortion of torus knots is unbounded,
${\delta}(T_{p,q})>\frac{1}{160} \, {\rm min}(p,q)$,
answering a 1983 question of M. Gromov \cite{Gromov}. A different construction of families of knots with unbounded 
geometry using the volume of branched covers is due to Gromov-Guth \cite{GG}.

The {\em intrinsic distortion} of a complex $K$ with respect to embeddings in ${\mathbb R}^d$ is defined as
\begin{equation} \label{intrinsic distortion}
D_d(K)\; :=\; \inf_i \; {\delta}(i(K)),
\end{equation}
where the infimum is taken over all embeddings ${i:K \hookrightarrow {\mathbb R}^d}$.
It is established in \cite{Bourgain} that any $N$-point metric space can be embedded in the Euclidean space
of dimension $O(log\, N)$ with distortion $O(log\, N)$. The vertex sets of expander graphs (discussed in more detail in section
\ref{distortion section}) of bounded degree with the graph metric are examples of metric spaces
for which this distortion bound is tight (even for embeddings into the Hilbert space $l^2$) \cite{LLR}. 

However if one is interested in embeddings into the Euclidean space of a {\em fixed} dimension $d$, the bounds on distortion are polynomial: 
\cite{Matousek} gives upper and lower bounds on distortion of an $N$-point metric space into ${\mathbb R}^d$, $d\geq 3$ (recall that $f(N)={\Omega}(g(N))$ means: $f(N)>c\, g(N)$ asymptotically):
\begin{equation} \label{Oomega bounds}
O(N^{2/d}log^{3/2}N)\; \, {\rm and }\; \, {\Omega}(N^{2/d}) \; (d \; \, {\rm even}),\; \,  {\Omega}(N^{\frac{2}{d+1}}) \; \, (d\; \, {\rm odd}).
\end{equation}

In section \ref{distortion section} we analyze the distortion of embeddings of expander graphs in a Euclidean space of a fixed dimension, where
the metric space is taken to be the entire graph (including its edges), not just the vertex set.
We give a geometric argument, to some extent
following the natural approach of \cite{KB, GG} of slicing the embedding by codimension $1$ hypersurfaces, to show that for families of expander graphs of bounded degree the
distortion $D_d$ with respect to embeddings into ${\mathbb R}^d$ has a lower bound ${\Omega}(N^{1/(d-1)})$, see theorem \ref{expander theorem}. 
It is an interesting question whether there is a transition
to exponential distortion for embeddings of complexes below the stable range, similar to the behavior of thickness  discussed in section \ref{n-complexes};
for intrinsic distortion we do not have a candidate family to propose.

\section{Non-recursive complexity of embeddings in codimensions $0,1,2$.} \label{relative section}

It is proved in \cite{MTW} that the embedding problem  for both $d$- and $(d-1)$- complexes in ${\mathbb R}^d$ is undecidable, for $d\geq 5$. Their starting point is Novikov's theorem that there  can be no decision procedure for determining if a $d$-complex ${\Sigma}^d$, $d\geq 5$, which we may assume is a homology $d$-sphere, is actually homeomorphic to $S^d$. A second ingredient is Newman's theorem \cite{Newman} that any PL embedding (in the technical sense: a piecewise-linear embedding of a subdivision of the standard triangulation of) $S^{d-1}\hookrightarrow S^d$ is topologically bicollared. Then the topological Schoenflies theorem implies that the closed complementary domains of any such embedding are homeomorphic to $B^d$, the $d$-ball. From this it follows directly that there can be no algorithm  either for determining the embeddability of ${\Sigma}^d\smallsetminus {\rm int}({\Delta}^d)$  or ${\Sigma}^d\smallsetminus \sqcup_i {\rm int}({\Delta}_i^d)$ in ${\mathbb R}^d$, where the first complex is a homology  $d$-ball obtained by deleting the interior of a single $d$-simplex and the second complex is the $(d-1)$-skeleton ${\rm Sk}^{d-1}({\Sigma}^d)$. Let 
$\{ J_i^d\}=\{ {\Sigma}_i^d\smallsetminus {\rm int}({\Delta}^d)\}$ 
and
$\{ J_i^{d-1}\}=\{ { \rm Sk}^{d-1}({\Sigma}_i^d)\}$ denote these undecidable families. 

The following theorem, giving a refinement complexity (rc) - analogue of these statements, is given in \cite[Corollary 1.2]{MTW}. Modulo a technical conjecture \ref{conj} below,
the statement also holds for thickness, see conjecture \ref{conj1}.

\begin{theorem} \label{non-recursive thm} \sl
For each $d\geq 5$ there exists a sequence $\{ K_i\}$ of $d$-complexes ($(d-1)$-complexes) with $n_i$ simplices which after subdivision PL 
embed in ${\mathbb R}^d$, but there is no recursive function of $n_i$ which lower bounds ${\rm rc}(n_i)$.
\end{theorem}

{\em Proof.} The Tarski-Seidenberg theorem \cite{Tarski} on quantifier elimination implies that arbitrarily quantified sentences in first order semi-real algebraic geometry are decidable. Clearly the existence of a linear embedding (no refinement) of a fixed finite complex in ${\mathbb R}^d$ can be reduced to quantified polynomial inequalities. 
For example the statement that there are disjoint line segments in ${\mathbb R}^3$ may be expressed as:

$$\exists\;  \vec{x_1}, \vec{x_2},  \vec{x_3},  \vec{x_4},  \; \forall \; s, t\; \; F(\vec{x_1}, \vec{x_2},  \vec{x_3},  \vec{x_4}, s, t)>0,\;
0\leq s \leq 1,\; 0\leq t\leq 1,\;\, {\rm where} $$
$$
F(\vec{x_1}, \vec{x_2},  \vec{x_3},  \vec{x_4}, s, t)\, =\, 
\parallel s \vec x_1+(1-s) \vec x_2-t \vec x_3-(1-t) \vec x_4\parallel^2.$$

Hence ``linear embedding'' is decidable. Let $\{ J_{i_j} \}$ be the subsequence of either family of complexes (dimension $d$ or $d-1$) discussed prior to theorem \ref{non-recursive thm} which actually {\em do} have PL embeddings into ${\mathbb R}^d$. Consider the function ${\rm rc}(n_{i_j})$. If this function were upper bounded on $\{ n_{i_j} \}$ by any recursive function 
$r(n_i)$ we could write an algorithm  for deciding the PL embeddability of the family $\{ J_i\}$. We would run the (in principle) Tarski-Seidenberg 
algorithm on each of the finitely many subdivisions  of $J_i$ containing $\leq m$ simplices for each $m\leq r(n_i)$. If one of these 
subroutines found a solution to the quantified embedding inequalities we would know that $J_i$ PL embed in ${\mathbb R}^d$. If by the time we have checked all subdivisions with $\leq r(n_i)$ simplices  no solution has been found we could halt the search knowing further subdivision to be fruitless. $\{ J_{i_j} \}$ serves as the $\{ K_j\}$ in the statement.
This contradiction proves the theorem. 
\qed

\begin{conjecture} \label{conj1} \sl For each $d\geq 5$ there exists a sequence $\{ K_i\}$ of $d$-complexes ($(d-1)$-complexes) with $m_i$ simplices which can be piecewise smoothly embedded in ${\mathbb R}^d$ but the intrinsic thickness satisfies: $${\rm thickness}\, (K_i)\; < \; r(n_i)$$
where $r$ is any positive recursive function.
\end{conjecture}

To discuss the relation between theorem \ref{non-recursive thm} and conjecture \ref{conj1}, note that
our definition \ref{defi thickness} of piecewise smooth embeddings is rather minimal.
It has the virtue that because the closed simplices are smoothly embedded, the complex automatically satisfies the Whitney stratification conditions (A and B).
However in other regards our definition is bare bones. It does not control how
closely adjacent simplices may approach each other.
For example our definition allows the $1$-complex with two $1$-simplices joined to form the letter V to be embedded in ${\mathbb R}^2$ by:
$$f_1(s)=(s,0), \; 0\leq s\leq 1,$$ $$f_2(t)=(t,F(t)), \; 0\leq t\leq 1,$$

where $F(0)=0$, $F(1)=1$ is any smooth function ($F(t)>0$ for $t>0$). $F$ might satisfy $F(1/n)=({\rm Ackermann}(n))^{-1}$, or even decay nonrecursively.

While the minimal nature of our definition  is adequate for the main results in section ~\ref{n-complexes}, it makes it quite technical (but we think possible) to prove a lemma
relevant to theorem \ref{non-recursive thm} in the context of thickness. We choose to state this lemma as conjecture \ref{conj} which we recommend to anyone with technical interest in stratified sets:

\begin{conjecture} \label{conj}
If $K^n$ admits a piecewise smooth embedding of thickness $\epsilon$ into ${\mathbb R}^d$ (more specifically,
into the unit ball $B^d_1\subset {\mathbb R}^d$) then ${\rm rc}(K)< r({\epsilon})$ where $r$ is a recursively computable 
function of ${\epsilon}$.
\end{conjecture}

{\bf Remark}. Given this, conjecture \ref{conj1} follows from the proof of theorem \ref{non-recursive thm}.

Next we establish non-recursive complexity in codimension 2 for a family of relative embedding problems.

\begin{theorem} \label{relative theorem} \sl
For each $d\geq 5$ there exists a sequence of $(d-2)$-complex pairs 
$(K_i, L_i)$ with a fixed embedding $L_i\subset {\mathbb R}^d$
such that the  refinement complexity (and conjecturally the thickness) of any extension to an embedding $K_i\subset {\mathbb R}^d$
is a non-recursive function of the number of simplices of $K_i$.
\end{theorem}

{\em Proof.}
Indeed, consider a $2$-complex $K$ realizing a finite presentation of a group with undecidable word problem. Some thickening of $K$ to a  
$d$-dimensional $2$-handlebody $H$, $d\geq 5$,  embeds into the $d$-sphere $S^d$.
There is a sequence of curves $C_i$ in $H$, say those exhausting the kernel
$[{\pi}_1(1$-skeleton of $K)\longrightarrow {\pi}_1(K)]$, bounding disks which 
embed in $H$, so that the refinement complexity of such embeddings is necessarily non-recursive as a function of the length of the boundary curve. Otherwise an algorithm  which searched by brute force for such disks  up to a (recursively) computed upper bound on complexity would succeed in solving the word problem in ${\pi}_1(K)$.

Consider the complement $X:= S^d\smallsetminus H$, a codimension zero subcomplex
of $S^d$.
For each $i$ the curve $C_i$ may be assumed to be embedded in the $1$-skeleton of
$\partial X$ for a sufficiently fine triangulation.
Define $K'_i=X\cup_{C_i} D^2$. Then $(K'_i, X)$ is a $d$-dimensional pair satisfying the conclusion of the theorem.

The dimension of the complex can be reduced to $d-2$ by considering $Y:=(d-1)$-skeleton of $X$,
with the same fixed embedding into $S^d$.
Repeat the construction above with $Y$ in place of $X$, yielding pairs $(K_i, Y)$. 
Now the group in which we have to study the word problem is the original group 
${\pi}_1(K)$ free product with a finitely generated free group. The free summand does not affect the decision problem.
\qed 

The  codimension $2$ relative examples $(K_i, L_i)$ with a fixed embedding of $L_i$ in theorem \ref {relative theorem} do not admit an immediate generalization
to the absolute setting  $L_i=\emptyset$. Indeed, fixing the fundamental group is 
a crucial ingredient in the proof, and if one considers arbitrary embeddings of a complex into ${\mathbb R}^d$ the fundamental group of the complement  is not controlled by Alexander duality. (The Stallings theorem gives a partial control of the fundamental group, modulo
its lower central series for some codimension $2$ embeddings, see section \ref{2-complexes} for more details.)

\section{An upper bound for $n$-complexes in ${\mathbb R}^{2n}$.} \label{rc bound}

By general position a PL embedding of an $n$-complex $K$ in the stable range does not require any subdivisions.
At the other extreme discussed in section \ref{relative section}, in codimensions $0,1$ (and in codimension $2$ in the relative setting), 
the refinement complexity (definition \ref{rc}) is non-recursive. In this section we 
establish an upper bound at the edge of the stable range, for $n$-complexes in ${\mathbb R}^{2n}$:

{\bf Theorem \ref{rc upper bound thm}.} {\sl
Let $K$ be an $n$-complex with $N$ simplices which topologically embeds in ${\mathbb R}^{2n}$, $n\geq 3$.
Then ${\rm rc}(K, 2n)=O(e^{N^{4+{\epsilon}}})$ for any ${\epsilon}>0$.}

To be more specific, the upper bound will be obtained in terms of $N=$max(number of $(n-1)$-simplices, number of $n$-simplices).
Examples in section \ref{n-complexes} complement this upper bound on refinement complexity with an exponential lower bound.

{\em Proof of theorem \ref{rc upper bound thm}.} We start by briefly recalling the definition of the van Kampen obstruction and refer the reader to \cite{FKT} for more details. Predating a formal definition of cohomology,  
Van Kampen's 1933 paper \cite{van Kampen} gave a rough description of a cohomological embedding obstruction 
for $n$-complexes $K^n$  in ${\mathbb R}^{2n}$. Any general position map $f\co K^n\longrightarrow {\mathbb R}^{2n}$
gives rise to a ${\mathbb Z}/2$-equivariant $2n$-dimensional simplicial cochain $o_f$ on $K\times K\smallsetminus {\Delta}$, where $\Delta$ is the union
of all products of simplices that have at least one vertex in common. Given $2n$-dimensional ${\sigma}\times {\tau}\in K\times K\smallsetminus{\Delta}$, 
$o_f({\sigma}\times{\tau})$ is defined as the algebraic intersection number $f({\sigma})\cdot f({\tau})$ in ${\mathbb R}^{2n}$.
The cohomology class $o(K)\in H_{{\mathbb Z}/2}^{2n}(K\times K\smallsetminus{\Delta};{\mathbb Z})$ of $o_f$ is independent of the map $f$: 
a general position  homotopy between two maps
$f, f'$ has finitely many non-generic times $\{ t_i\}$ when an $n$-cell intersects an $(n-1)$-cell, and it is not difficult to see
that $o_f$, $o_{f'}$ differ by the coboundary of the sum over $\{ t_i\}$ of the products
$(n$-cell$)_{t_i} \times ((n-1)$-cell$)_{t_i}$
 \cite{FKT}. (Note that the action of ${\mathbb Z}/2$ on the cohomology coefficients is trivial \cite{Shapiro, Melikhov} and not $(-1)^n$ as stated in \cite{FKT}.)
Moreover, for any cochain representative $c$ of $o(K)$ there is a map $f$ such that $c=o_f$ \cite[Section 2.4, Lemma 2]{FKT}.

The vanishing of the van Kampen obstruction is clearly a necessary condition for embeddability of $K^n$ in ${\mathbb R}^{2n}$ for any $n$.
For $n>2$ it is also a sufficient condition, due to the validity of the Whitney trick, see \cite{Shapiro}, \cite[Theorem 3]{FKT} and the
discussion below. (For $n=2$ the obstruction is incomplete \cite{FKT}.)

Suppose an $n$-complex $K$ with $N$ simplices embeds in ${\mathbb R}^{2n}$. 
First map it linearly into ${\mathbb R}^{2n}$, in general position. That is, consider a generic map
of the vertices of $K$ to ${\mathbb R}^{2n}$ and extend linearly to get a map $f\co K\longrightarrow {\mathbb R}^{2n}$. 
The rank of the $2n$-th cellular cochain group $C^{2n}$ of $K\times K\smallsetminus {\Delta}$ is $<N^2$, and the possible values of the van Kampen cochain $o_f$ 
on each $2n$-cell are $0, \pm 1$, so $|| o_f ||< N$. Here and in other estimates below we will work with the $l^2$-norms on the cochain groups
$C^{2n}$, $C^{2n-1}$ where the $2n$- (respectively $(2n-1)$-) dimensional product cells ${\sigma}\times {\tau}$ are taken as orthonormal generators.

There are two steps in changing this map to a PL embedding. Since $K$ embeds in ${\mathbb R}^{2n}$, the van Kampen obstruction $o(K)$ vanishes and 
$o_f$ is a coboundary. Therefore, as discussed above, there are ``finger move'' homotopies pushing $n$-simplices of $K$ across $(n-1)$-simplices, 
so that the result is a PL map $g\co K\longrightarrow {\mathbb R}^{2n}$ with $o_{g}\equiv 0$, in other words
$g({\sigma})\cdot g({\tau})=0$ for any two non-adjacent $n$-simplices ${\sigma}, {\tau}$ of $K$. 
Each finger move homotopy amounts to taking a connected sum of 
an $n$-simplex $f({\sigma})$ with a small oriented PL $n$-sphere linking $\pm 1$ an $(n-1)$-simplex $f(\nu)$ of $K$ (cf. figure 2.1 in \cite{FKT}). 
The connected sum is taken along a straight line segment in general position with $f(K)$.
This move requires a constant (depending on the dimension $n$) number of subdivisions of the simplex $\sigma$. 

The second step consists of applications of the Whitney move to remove the intersections among the
images of the $n$-simplices under $g$. In general there are three types of Whitney moves: corresponding to 

(1) a $\pm$ pair of intersections of two non-adjacent
simplices $g({\sigma}^n), g({\tau}^n)$,  

(2) an intersection point of two distinct adjacent $n$-simplices, and 

(3) a self-intersection point of an $n$-simplex.

To begin with, the map $f$ is linear on each simplex and so there are no intersections of types (2), (3).
The finger moves implementing a homotopy to $g$ may introduce intersections of type (2), however since the finger move
always involves non-adjacent simplices ${\sigma}^n, {\nu}^{n-1}$, self-intersections will not be introduced and so the move
of type (3) is not needed.

 Correspondingly, there are two types of Whitney circles $C$:  the usual Whitney circle consisting of two arcs pairing up the
intersection points in two non-adjacent simplices ${\sigma}^n, {\tau}^n$, and a slightly different type   connecting an intersection between two
adjacent simplices through a common vertex. In either case, the Whitney circle consists of a bounded number 
(to be specific, at most $6$) line segments.
For the Whitney disk $W$ in each case take the cone on the
Whitney circle $C$. A generic choice of the cone point ensures that the interior of $W$ is disjoint from $g(K)$. Therefore 
each Whitney disk $W$ consists of at most 6 $2$-simplices. 

Due to the bounded combinatorial complexity of the Whitney disk discussed above,
the classical Whitney move of type (1), changing one of the simplices $g({\sigma})$, $g({\tau})$ by an ambient PL isotopy, requires
an  a priori bounded number of subdivisions.
The Whitney move of type (2) is described in lemma 5 in \cite{FKT}. It eliminates an intersection point between two adjacent simplices while
also introducing a bounded number of subdivisions. 

Since each finger move gives rise to at most $N$ new intersection points, the number of Whitney moves
is bounded by $N$ times the number of finger moves.
The remaining problem is to get an estimate on the number of finger moves 
needed to get from $f$ to $g$. This problem is equivalent to getting an upper bound on the $l^1$-norm of a 
$(2n-1)$-cochain which maps onto the van Kampen $2n$-cocycle $o_f$. We will work with the $l^2$-norms below,
however the required estimate on the $l^1$-norm will follow since $||v||_1\leq N ||v||_2$ (recall that the
ranks of $C^{2n-1}, C^{2n}$ are $<N^2$) . 

Denote by $A$ the matrix representing the coboundary map $C^{2n-1}\longrightarrow C^{2n}$ with 
respect to the chosen cell basis. The minimal norm of a cochain $c\in C^{2n-1}$ such that $Ac=o_g$ is upper bounded by 
${\lambda}^{-1}||o_f||$, where
\begin{equation} \label{lambda}
{\lambda}={\min}_{v\neq 0\in ({\rm ker}\, A)^{\perp}\subset C^{2n-1}}\; \frac{||Av||}{||v||}\; =\; {\min} \; \frac{\langle A^t Av , v\rangle^{1/2}}{||v||}.
\end{equation}

Therefore we need to estimate the smallest absolute value ${\mu}$ of a non-zero eigenvalue of the symmetric matrix $M:= A^t A$.
The coboundary map from $(n-1)$- to $n$-cochains of the original complex $K$
has $O(N)$ $\pm 1$ entries, the rest are zeros. The matrix $A$ has $O(N^2)$ $\pm 1$'s, the rest of the matrix entries are zeros. 
Then the sum of the absolute values of the entries of $M$ is bounded in absolute value by $O(N^4)$.

Let $\chi$ be the characteristic polynomial of $M$ with variable $x$. Unless the original linear map happened to be embedded and no subdivision was necessary, 
$M$ must have a non-zero eigenvalue, otherwise the non-zero obstruction element could not be a coboundary. It follows that $\chi$ must have at least two non-zero 
terms; let the term of lowest degree read:   $c_k x^k$. To estimate how small (in absolute value) a root of $\chi$ could possibly be, we take $c_k$ as small in absolute 
value as possible, $c_k = -1$ and all the coefficients of higher terms as large as possible. Estimating this bound $B$ by computing determinants we find the maximum 
is  $c_j = O((N^4)!N^4),\,  N^2>j>k$.  We get a (smallest case) equation:
$x^k = \sum_{i=k+1}^{N^2}  B x^i$.

Summing the geometric series  obtain $1 = Bx/(1-x)$.  
We find a lower bound on $\mu$:
$\mu$ can be no smaller than $(O((N^4)!N^4))^{-1}$.
Then the minimal norm of a cochain $c$ with $Ac=o_g$ is 
$<{\lambda}^{-1}||o_f||<{\mu}^{-1/2}N=O(e^{N^{4+{\epsilon}}}).$
\qed

{\bf Remarks}. 1. The upper bound on ${\rm rc}(K, 2n)$ in theorem \ref{rc upper bound thm} is established in terms of the weakest possible
isoperimetric constant of the coboundary map $C^{2n-1}\longrightarrow C^{2n}$ of the deleted product $K\times K\smallsetminus {\Delta}$ that works for any $K$.
If one is interested in a specific $n$-complex $K$ with a larger isoperimetric constant, the proof above would give 
a more optimal bound on the refinement complexity of $K$.

2. A more careful linear algebra estimate may improve the power on $N$ in the exponent in the statement of the theorem. 
It does not seem likely however that a minor improvement along these lines would close the gap with the lower bound $c^N$ proved in the
next section.

3. It seems reasonable to believe that a version of theorem \ref{rc upper bound thm} also holds for thickness. Such a statement would follow if one
established an estimate on thickness in terms of the number $M$ of simplices in a PL embedding. For example, it is plausible that 
thickness $\sim M^{-1}$.

\section{Exponentially thin $n$-complexes in ${\mathbb R}^{2n}$.} \label{n-complexes}

This section gives the proof of theorem \ref{exponential theorem}.
Note that most $n$-dimensional complexes do not embed
into ${\mathbb R}^{2n}$, and for those complexes that do embed a given generic map
into ${\mathbb R}^{2n}$ may not in general  be approximable by embeddings.
The idea of the proof is to construct a sequence of $n$-complexes  which
``barely'' embed into ${\mathbb R}^{2n}$: any embedding
necessarily involves an exponential amount of linking, and as a consequence
an exponential bound on thickness and recursive complexity.

Let $K^n_0$ denote 
the $n$-skeleton of the $(2n+2)$-simplex, $({\Delta}^{2n+2})_n$, 
with a single $n$-simplex $T$ removed. 
Denote the boundary 
$(n-1)$-sphere of the
missing simplex by $S^{n-1}_1$, and let $S^n_2$ be the $n$-sphere 
spanned in $K$ by the $n+2$ vertices which are not in $S_1$. 

\begin{proposition} \label{mod 2 linking}  \sl
The $n$-complex $K^n_0$ embeds into ${\mathbb R}^{2n}$. Moreover, for any embedding 
$i\! :K^n_0\longrightarrow {\mathbb R}^{2n}$ the mod-$2$ linking number ${\rm lk}_{{\rm mod\,  2}}(i(S^{n-1}_1), i(S^n_2))$ 
 is non-zero.
\end{proposition}

This proposition essentially follows from the work of van Kampen \cite{van Kampen}. We give an outline of the argument below and refer the reader to lemma 6 in \cite{FKT} for more details.  (\cite{FKT} states the proof for $2$-complexes in ${\mathbb R}^4$; the argument for all $n\geq 2$ is directly analogous.)

The complex $({\Delta}^{2n+2})_n$ for $n=1$  is the complete graph on $5$ vertices (well-known from the Kuratowski planarity criterion), which can be drawn on the plane with a single intersection point between two non-adjacent edges.
Van Kampen \cite{van Kampen} (and independently Flores \cite{Flores}) observed that this is also true for $n\geq 2$: $({\Delta}^{2n+2})_n$ can be mapped into ${\mathbb R}^{2n}$ with a single double point between two non-adjacent $n$-simplices. For example, for $n=2$ this can be seen concretely as follows: consider the $2$-skeleton ($\cong S^2$)
of ${\Delta}^3\subset {\mathbb R}^3\times 0 \subset {\mathbb R}^4$. To get the $2$-skeleton
of ${\Delta}^6$ add another vertex $v_5$ inside the $2$-sphere $\partial {\Delta}^3$ in ${\mathbb R}^3$, and $v_6, v_7$ in ${\mathbb R}^4_+,  {\mathbb R}^4_-$ respectively.
Connect $v_6, v_7$ by an edge intersecting ${\mathbb R}^3\times 0$ in a point outside
of ${\Delta}^3$. It is not difficult to see that all $2$-cells embed disjointly, except for
the $2$-cell with vertices $v_5 v_6 v_7$ whose boundary circle links the  $2$-sphere
$\partial {\Delta}^3$. Adding this $2$-cell introduces a single double point with one of the four faces forming the $2$-sphere.

Van Kampen proved that the embedding obstruction (discussed in the proof of theorem \ref{rc upper bound thm} in section \ref{rc bound}) is non-trivial for the complex  $({\Delta}^{2n+2})_n$
by starting with the map with a single double point discussed above and showing that any homotopy preserves (mod 2) the sum of the
 cochain $o_f$ over all $2n$-cells.
Now suppose the linking number is zero in the statement of the proposition for some embedding $i\co K_0\hookrightarrow {\mathbb R}^{2n}$.
Attach an $n$-cell $T$ to $S^{n-1}_1$ and map it into ${\mathbb R}^{2n}$, thus extending $i$ to a map
$f\co ({\Delta}^{2n+2})_n\longrightarrow {\mathbb R}^{2n}$.
 By assumption the intersection number
$f(T)\cdot f(S^n_2)={\rm lk}(i(S^{n-1}_1), i(S^n_2))$ equals  zero (mod 2). 
Then the sum of the values of the cochain $o_f$ over all $2n$-cells is zero, a contradiction.
This concludes an outline of the proof of proposition \ref{mod 2 linking}.
\qed

We are now in a position to construct the complexes used in the proof of
theorem \ref{exponential theorem}. For convenience of the reader we restate it here:

\smallskip

{\bf Theorem \ref{exponential theorem}.} {\sl
For each $n\geq 2$ there exist families of $n$-complexes $\{ K_l\} $ with
$m_l\longrightarrow \infty$ simplices and bounded local combinatorial complexity which embed into ${\mathbb R}^{2n}$ and the thickness of any such embedding is at most $c^{-m_l}$. 
Moreover,
\begin{equation} \label{rc expo}
 C^{m_l}\; <\; {\rm rc}(K_l, 2n)\; < \; \infty.
\end{equation}
Here the constants $c, C>1$ depend on $n$.}

\smallskip

{\em Proof.} Fix $n$ and for each $l\geq 1$ define the $n$-complex $K_l$ to be the $l$-fold mapping telescope 
\begin{equation} \label{telescope} S^{n-1}\overset{\times 2} \longrightarrow \, S^{n-1} \overset{\times 2} \longrightarrow \, \ldots \overset{\times 2}\longrightarrow  \, S^{n-1} \overset{\times 2}\longrightarrow K_0
\end{equation}

where the last arrow denotes the degree $2$ map from $S^{n-1}$ to the $(n-1)$-sphere $S^{n-1}_1\subset K_0$ from the statement of proposition \ref{mod 2 linking}. 
Sometimes the dimension of this sphere will be dropped and it will be denoted $S_1$.

Concretely each stage of the mapping telescope is a copy of $S^{n-1}\times[0,1]$
triangulated so that $S^{n-1}\times 0$ receives an induced triangulation as
$\partial {\Delta}^n$, and $S^{n-1}\times 1$ is triangulated as the $2$-fold branched cover of $\partial {\Delta}^n$ along $\partial {\Delta}^{n-2}$, for some ${\Delta}^{n-2}\subset \partial{\Delta}^n$. The map between each product segment of the telescope is the aforementioned  branched cover.
Therefore $K_l$ has $O(l)$ $n$-simplices and the local combinatorial complexity of $K_l$ does not increase with $l$.

Next observe that the embedding $i\co K_0\hookrightarrow {\mathbb R}^{2n}$ discussed in proposition \ref{mod 2 linking} extends to an embedding $K_l\hookrightarrow {\mathbb R}^{2n}$.  Recall that the $n$-cell $T= ({\Delta}^{2n+2})_n\smallsetminus K_0$ bounded by $S_1$ can be mapped into ${\mathbb R}^{2n}$ creating a single double point
with a non-adjacent $n$-cell. In particular, a collar $S_1\times [0,1]$ on the boundary $(n-1)$-sphere of $T$ embeds into ${\mathbb R}^{2n}$ disjointly from the rest 
of the complex $K_0$. A disk normal bundle over $T$ in ${\mathbb R}^{2n}$, restricted to the collar,  which we will parametrize as $S_1\times D^n\times[0,1]$,
also embeds into ${\mathbb R}^{2n}$ disjointly from the rest of $K_0$. 
The sphere $S^{n-1}_1\subset K_0$ is identified with $S_1\times 0\times 0$ in this normal bundle.

By general position a degree $2$ map $S^{n-1}\longrightarrow
S^{n-1}_1$ can be perturbed to an embedding $S^{n-1}\hookrightarrow S^{n-1}_1\times D^n\times t_1$, 
$0<t_1<1$. Then the mapping cylinder $S^{n-1} \overset{\times 2}\longrightarrow S_1$ embeds
level-wise into $S_1\times D^n\times[0,t_1]$, where the target sphere is identified with 
$S_1\times 0\times 0$. Proceeding by induction and embedding the spheres $S^{n-1}$ 
in $S_1\times D^n\times t_i$, $0< t_1<\ldots < t_l\leq 1$, the $l$ cells forming the mapping telescope (\ref{telescope}) embed in $S_1\times D^n\times [0,1]$, showing that
$K_l$ embeds into ${\mathbb R}^{2n}$.

To establish the bound on embedding thickness of the complexes $K_l$ claimed in theorem \ref{exponential theorem}, consider any
simplex-wise smooth embedding into the unit ball in ${\mathbb R}^{2n}$, $i\co K_l\longrightarrow B_1^{2n}$. Denote by $\bar S^{n-1}\subset K_l$ the left-most $(n-1)$-sphere
in (\ref{telescope}). 
According to proposition \ref{mod 2 linking},
the mod-$2$ linking number ${\rm lk}_{{\rm mod\,  2}}(i(S^{n-1}_1), i(S^n_2))$ is non-zero. Switching to the integer-valued linking number, it follows that
\begin{equation} \label{expo linking}
 | {\rm Lk}\, ( i(\bar S^{n-1}), i(S^n_2)) | \,  \geq \, 2^l.
\end{equation}

Recall from the introduction that  the thickness of an embedding $i\co K^n \longrightarrow B^{2n}_1\subset {\mathbb R}^{2n}$  
is defined as the supremum of $T$ such that 
the distance between the images of any two non-adjacent 
simplices is at least $T$ and all simplices have embedded $T$-normal bundles.
Next we will use the Gauss linking integral to get an upper bound on $T$.

For a $2$-component link in $S^3$ whose components are parametrized by
${\alpha}, {\beta}\co S^1\longrightarrow {\mathbb R}^3$,
the classical Gauss integral in 3-space computes the
linking number $L$ in terms of the triple scalar product of $\dot{\alpha}$, $\dot{\beta}$ and $v:={\alpha}(s)-{\beta}(t)$: 
$${\rm Lk}\, ({\alpha}, {\beta})\; =\; \frac{1}{4{\pi}}\, \int \frac{[\dot{\alpha}(s), \dot{\beta}(t), v]}{|v|^3}\, ds\, dt $$

Denoting by ${\alpha}, {\beta}$ parametrizations  of the spheres $i(\bar S^{n-1}), i(S^n_2)$ in ${\mathbb R}^{2n}$, 
consider the higher-dimensional analogue (cf. \cite{SV} for a more detailed discussion) of the Gauss integral computing
the degree of the corresponding map $S^{n-1}\times S^n\longrightarrow S^{2n-1}$, equal (up to a sign) to the linking number: 
\begin{equation} \label{Gauss} 
{\rm Lk}\, ({\alpha}, {\beta}) \; =\; c\int_{S^{n-1}\times S^n} \frac{{\rm det}\, [\, d{\alpha}, \, d{\beta},\,  v\, ]}{|v|^{2n}}
\end{equation}

Since $i$ is an embedding into the unit ball $B^{2n}_1$ and the embedding thickness is $T$, it follows that 
$T<|v|<1$. Then the integral (\ref{Gauss}) has an upper bound in terms of the $(n-1)$-volume $V_{\alpha}$,
the $n$-volume $V_{\beta}$ and the thickness $T$, leading to the inequality
\begin{equation}
2^l\; \leq \; {\rm Lk}\, ({\alpha}, {\beta})  \; < \; \frac{V_{\alpha}\, V_{\beta}}{T^{2n}}.
\end{equation}

Recall that that the spheres $\bar S^{n-1}, S^n_2$ are considered with standard triangulations.
Then 
\begin{equation} \label{inequality equation}
T^{2n} \, 2^l \, < \, V_{\alpha} \, V_{\beta}\, <\, (n+1)(n+2) \, V'_{\alpha}\,  V'_{\beta},
\end{equation}

where the prime denotes the largest volume among the simplices forming the two linking spheres.

We will next obtain a lower bound on the $2n$-volume of the $T$-normal bundles over the relevant simplices  of $\bar S^{n-1}, S^n_2$  in terms of their
$(n-1)$-, respectively $n$-volumes $V'_{\alpha}, V'_{\beta}$. There is a polynomial expression, the  classical Weyl tube formula  \cite{Weyl}, for the volume of an $\epsilon$-regular neighborhood 
of a submanifold of a Euclidean space in terms of the volume of the submanifold.
For example, the Weyl formula for closed surfaces in ${\mathbb R}^3$ states:
$${\rm Vol} ({\mathcal N}_{\epsilon} {\Sigma})\; =\; 2\, {\rm Area} ({\Sigma})\, {\epsilon} \, +\, \frac{4{\pi}}{3}\,  {\chi}({\Sigma})\, {\epsilon}^3,$$

cf. \cite[section 1.2]{Gray}. Note that in general this formula
does not immediately give a lower bound for the volume in terms of the area since higher order
terms may be negative (and we don't have an apriori estimate on how small $\epsilon$ is).
Instead we will use a more direct argument to get a rough estimate relating 
the volumes of the submanifold and of its normal bundle.

For brevity of notation let $S$  denote the relevant $(n-1)$- (respectively $n$-) simplex in the sphere $i({\bar S}^{n-1})$, respectively in $i(S_2^n)$, and let $q$ denote its dimension.
Let $\nu$ be the distance from $S$ to its nearest focal point and set ${\epsilon}={\nu}/2$. 
Consider the restriction $\phi$ of the exponential map (defined on the tangent bundle of ${\mathbb R}^{2n}$) 
to the ${\epsilon}$-disk normal bundle $N_{\epsilon}$ over $S$.
Denote the ${\epsilon}$-regular neighborhood of $S$ in 
${\mathbb R}^{2n}$ by ${\mathcal N}_{\epsilon}$. 

\begin{figure}[ht]
\includegraphics[width=5cm]{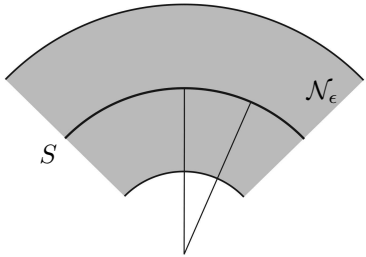}
{\Small
    \put(-26,61){${\mathcal N}_{\epsilon}$}
    \put(-127,37){$S$}}
\caption{The ${\epsilon}$-tube around $S$.}
\end{figure}

We will derive a rough bound $|{\rm det}\, (\, d_x \, {\phi})|\, > \, 1/2^q$ on the differential of the diffeomorphism  ${\phi}\co N_{\epsilon} \longrightarrow {\mathcal N}_{\epsilon}$ 
at any point $x\in N_{\epsilon}$ to get an estimate:
\begin{equation} \label{volume estimate}
{\rm Vol}_{2n}({\mathcal N}_{\epsilon})\; \geq \; \frac{{\rm Vol}_{2n-q}(D^{2n-q}_{\epsilon})\, {\rm Vol}_q(S)}{2^q},
\end{equation}

The bound on the differential is particularly easy to see in the case of a hypersurface, pictured in figure 1:
due to the choice of ${\epsilon}={\nu}/2$ where ${\nu}$ is the distance to the nearest focal point,
the tangent vectors to $S$ do not contract by more than a factor $1/2$ when mapped to any point within ${\mathcal N}_{\epsilon}$.
(In some directions the tangent vectors expand - however we are looking for just a very  rough estimate so we take
the ``worst case'' bound in every direction.) This gives a lower bound of $1/2$ on the eigenvalues of the differential
restricted to vectors tangent to $S$, and the vectors tangent to the normal disk are mapped isometrically by
$d_p \, {\phi}$.

Proving (\ref{volume estimate}) in general (for submanifolds of higher codimension) requires a bit more care, and we will 
refer to \cite{Gray} for details of the statements below. Fix $p\in S$
and let $u$ be a unit vector normal to $S$ at $p$.  Denote by $S_t$ the boundary of the tube of radius
$t$ around $S$. Then the principal curvatures ${\kappa}_i, i=1,\ldots, 2n-1$ of $S_t$ at the point $p+tu$
satisfy the differential equation ${\kappa}'_i(t)\, =\, ({\kappa}_i(t))^2$ (\cite[Corollary 3.5]{Gray} states this in a more general setting
where $R^{2n}$ is replaced with an arbitrary Riemannian manifold.)

The first $q$ of these principal curvatures of $S_t$ may be thought of as corresponding to the directions tangent to $S$,
and they are given by ${\kappa}_i(t)={\kappa}_i(0)\, (1-t\, {\kappa}_i(0))^{-1}$, $i=1,\ldots, q$.
The rest 
correspond to the directions normal to $S$, ${\kappa}_i(t)=-t^{-1}$, $i=q+1, \ldots, 2n-1$ (they are equal to the corresponding principal
curvatures of the $t$-sphere normal bundle, since the exponential map is an isometry in these directions). Then as expected 
the determinant of the differential of the exponential map may be expressed in terms of the first $q$ principal curvatures:
\begin{equation} \label{det inequality}
|{\rm det}\, (\, d_{p+tu} {\phi})|\, = \, \prod_{i=1}^q \, (1-t{\kappa}_i(0)),
\end{equation}

see \cite[(3.27)]{Gray}. The principal curvatures ${\kappa}_i(0), i=1,\ldots, q$ satisfy the inequality 
${\kappa}_i(0)\, {\nu}\leq 1$, where as above ${\nu}$ denotes the distance to the nearest focal point 
\cite[section 8.1, Lemma 8.9]{Gray}. We assumed that $0\leq t\leq {\epsilon}={\nu}/2$, so $1-t{\kappa}_i(0)> 1/2.$ 
According to (\ref{det inequality}) then $|{\rm det}\, (\, d_p {\phi})|>2^{-q}$, concluding the proof of
the estimate (\ref{volume estimate}).

We proceed with the proof of theorem  \ref{exponential theorem}.
The complex $K_l$ is embedded in the unit ball and the $T$-normal bundles to the simplices are assumed to be embedded. Since 
${\nu}$ is the distance to a closest focal point, $T<  {\nu}=2{\epsilon}$. If $T>{\epsilon}$ then 
${\rm Vol}_{2n}({\mathcal N}_T)\, > \, {\rm Vol}_{2n}({\mathcal N}_{\epsilon})$. If $T<{\epsilon}$ then the inequality (\ref{volume estimate}) holds
for $T$ in place of ${\epsilon}$. In either case, as a consequence of (\ref{volume estimate}) one has
\begin{equation} \label{useful volume inequality}
{\rm Vol}_{2n}({\mathcal N}_T)\, >  \,  \frac{c_{2n-q} {(T/2)}^{2n-q}\, {\rm Vol}_k(S)}{2^q} \, > \, 
\frac{{T}^{2n-q}\, {\rm Vol}_k(S)}{2^{2n}},
\end{equation}

where $c_{2n-q}$ is the constant in the expression for the volume of the ball $D^{2n-q}_{T/2}$.
The effect of this constant is minor on the final estimate (\ref{thickness estimate}) on thickness below, 
so it is omitted in the last term in (\ref{useful volume inequality}).
To get a very rough
upper bound on ${\rm Vol}({\mathcal N}_T)$ note
that the $T$-regular neighborhood of each simplex is embedded in the ball of radius $2$, so 
its $2n$-volume is less than $ (2{\pi})^n/n!$. 
Combine this with the estimate (\ref{useful volume inequality}) where $S$ is largest volume simplex
in $i(\bar S^{n-1}), i(S^n_2)$ as in (\ref{inequality equation}),  so ${\rm Vol}_k(S)=V'_{\alpha}$, respectively $V'_{\beta}$:
\begin{equation} \label{bundle estimate}
\frac{T^{n+1}}{2^{2n}} \; V'_{\alpha} \, < \, \frac{{\pi}^n}{n!}, \;\; \; \frac{T^{n}}{2^{2n}} \; V'_{\beta} \, < \, \frac{{\pi}^n}{n!}
\end{equation}  

Combine (\ref{inequality equation}), (\ref{bundle estimate}):
\begin{equation}
T^{2n} \, 2^l \, < \, (n+1) (n+2) \; \frac{(4 {\pi})^{2n}}{T^{2n+1} \, (n!)^2}
\end{equation}

It follows that
\begin{equation}
T^{4n+1}\, < \, 2^{-l} (n+1)(n+2)\;  \frac{(4{\pi})^{2n}}{(n!)^2}
\end{equation}

To a good approximation this implies
\begin{equation} \label{thickness estimate}
T\,  < \, 2^{-l/(4n+1)}\, n^{-1/2}.
\end{equation}

This is the sought exponential thinness in terms of $l$.


{\em Proof of the bound (\ref{rc expo}) on refinement complexity}. Instead of computing linking($(k-1)$-sphere, $k$-sphere) by an integral formula, we can instead use 
intersection numbers in a generic projection to ${\mathbb R}^{2k-1}$. Since for two linearly embedded simplices of dual dimensions the possible intersection numbers are: $-1, 0$,
or $1$, there must be exponentially many pairs of simplices to account for the exponential linking number. 
Denote by $M_i$ the number of simplices in a PL embedding. Then $M_i^2> c^{m_i}\Rightarrow M_i>(\sqrt c)^{m_i}$. \qed

{\bf Addendum}. {\sl The statement of theorem \ref{exponential theorem} (both for thickness and rc) applies  
to $(n+k)$-complexes in ${\mathbb R}^{2n+k}$, for any $k\geq 0$.}

{\em Proof of addendum.} Consider the $n$-complexes $K_i$ and their embeddings in the $2n$-sphere $K_i\subset S^{2n}$  constructed in the proof of theorem \ref{exponential theorem}.
Applying (iterated) suspension, one gets embeddings ${\Sigma}^k K_i \subset {\Sigma}^k S^{2n}\cong S^{2n+k}$. Recall the exponential bound (\ref{expo linking}) 
on the linking number of the $(n-1)$-sphere $i(\overline S)$ and the $n$-sphere $i(S_2)$ 
for any embedding $i\! : K_i\subset S^{2n}$. Then for any embedding  $i^k\! :{\Sigma}^k K_i \subset S^{2n+k}$ one has the same exponential bound on the linking number of
$i^k( \overline S)$ and $i^k({\Sigma}^k S_2)$. 
By induction assume exponential linking is established between $i^j(\overline S)$ and $i^j({\Sigma}^j S_2)$ in $S^{2n+j}$. 
Now look at any embedding $i^{j+1}$ of ${\Sigma}^{j+1} K_i$ intersect a small ${\epsilon}$-radius sphere $S_{\epsilon}^{2n+j}$ about one of the embedded suspension points and apply the inductive hypothesis. This linking number in $S_{\epsilon}^{2n+j}$ must of course agree with ${\rm lk} (i^{j+1}(\overline S), i^{j+1}({\Sigma}^{j+1} S_2))$ in $S^{2n+j+1}$.
The rest of the proof is identical to that of theorem \ref{exponential theorem}.
\qed

Our results on thickness and refinement complexity are summarized in Table 1.

\begin{table}[h]
\begin{tabular}{ c | c c c c c c c c c c c c c c c }
 n$\diagdown$d  & 1 & 2 & 3 & 4 & 5 & 6 & 7 & 8 & 9 & 10 & 11 & 12 & 13 & 14 & 15 \\
\hline
1 & 1P & 1P & 1P & 1P & 1P & 1P &&&&&&&&& \\
2 & *   & 1P & rr   & EE & 1P & 1P &&&&&& {\large 1P} &&& \\
3 &     &   *  & rr   & EE & EE & ${\mathbb E}$E & 1P & 1P & 1P &&& \\
4 &     &      &   *  & EE & Nn & EE & EE & ${\mathbb E}$E & 1P & 1P & \\
5 &     &      &      &  *   & Nn & Nn & EE  & EE & EE & ${\mathbb E}$E & 1P & 1P & \\
6 &    &   {\large *}   &      &      & *    & Nn & Nn  & EE & EE & EE & EE & ${\mathbb E}$E & 1P & 1P & \\
7 &    &       &      &      &      & *    & Nn & Nn  & EE & EE & EE & EE & EE & ${\mathbb E}$E & 1P \\
\end{tabular}
\bigskip
\bigskip
\caption{The $(n,d)$ entry in the table is a pair (Refinement complexity, thickness) for $n$-complexes in
${\mathbb R}^d$. 
\medskip
Explanation of the notation: \newline
Refinement complexity: 1=no refinement, E=(at least) exponential refinement,
${\mathbb E}$=exponential lower bound and $O(e^{N^{4+{\epsilon}}})$ upper bound for $n$-complexes in ${\mathbb R}^{2n}, n\geq 3$, 
N=non-recursive refinement, r=we conjecture 
\medskip
to be recursive.\newline
Thickness: P=polynomial, E=(at least) exponential, n=non-recursive (modulo conjecture \ref{conj}), r=we conjecture to be recursive.
\medskip
\newline
$*=$ cannot embed since $d<n$.}
\end{table}

{\bf Remarks}. We conjecture decidable (recursive) behavior in $d=3$ in line with numerous recent results
in $3$-manifold topology. Cells (3,4) and (4,4) are filled in ``EE'' since one may trivially add an additional $3$- or
$4$-simplex to the examples in (2,4). Similarly all ``Nn'' are certainly at least ``NE''.

The interested reader may find it instructive to compare this table with 
table 1 on page 4 in \cite{MTW} which lists the algorithmic complexity of the embedding problem for simplicial $n$-complexes into ${\mathbb R}^d$.
In particular, it is interesting to note that the algorithmic complexity of the embedding problem for $K^n$ in ${\mathbb R}^{2n}$, $n\geq 3$ is polynomial (this amounts to checking 
whether van Kampen's cohomological obstruction vanishes \cite{MTW}), while we have shown that the geometric complexity (thickness and also
refinement complexity) of embeddings in these dimensions is at least exponential.

\section{Further examples of $2$-complexes in ${\mathbb R}^4$, and links in ${\mathbb R}^3$.} \label{2-complexes}

Section \ref{n-complexes} introduced a family of $n$-complexes  which ``barely'' embed into ${\mathbb R}^{2n}$:
any embedding is necessarily exponentially thin.
Here we examine the case $n=2$ in more detail. The main purpose of this section is to discuss and formulate a number of open questions 
about related families of $2$-complexes 
in ${\mathbb R}^4$ and also about classical links in $3$-space. Their definition is based on the lower central series
 and other more general iterated word constructions.

This discussion of $2$-complexes is motivated by two problems. First, the bound on thickness in section 
\ref{n-complexes} is based on the integral formula for the linking number. It is an interesting question whether similar bounds on thickness
may be obtained when the linking numbers vanish but a (suitably interpreted) higher-order Massey product is non-trivial. 
There are integral formulas that may be used to compute Massey products, but they do not seem to be directly suitable for
getting a thickness estimate as above. At the end of the section we also pose a related question for links in $3$-space where Massey products 
correspond to Milnor's $\bar\mu$-invariants \cite{Milnor}.

Another motivation is suggested by Milnor's representation \cite{Milnor} of a link which for convenience of the reader is 
reproduced in figure 2. This is a Brunnian link, in particular the sublink formed by the components $l_1,\ldots, l_{q+1}$ is the
$(q+1)$-component unlink, so the fundamental group of its complement 
is the free group generated by meridians $x_1,\ldots, x_{q+1}$ to these components.
The remaining curve $w'_q$ represents a $(q+1)$-fold commutator 
$w'_q\, =\, [x_1,[x_2,[x_3,\ldots, [x_q,x_{q+1}]]]]$ in this free group. 

\begin{figure}[ht]
\medskip
\includegraphics[height=2.5cm]{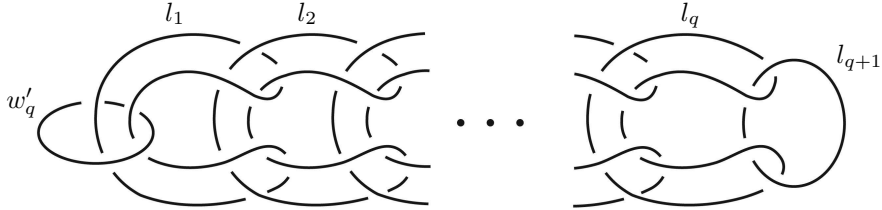}
{\Small
\put(-320,40){$w'_q$}
\put(-260,73){$l_1$}
\put(-210,73){$l_2$}
    \put(-65,73){$l_q$}
    \put(-6,57){$l_{q+1}$}}
\bigskip
\caption{$w'_q\, =\, [x_1,[x_2,[x_3,\ldots, [x_q,x_{q+1}]]]]$, where $x_i$ denotes a meridian to the corresponding link component $l_i$.}
\end{figure}

The length of $w'_q$
with respect to the word metric in the free group $F_{x_1,\ldots, x_{q+1}}$  is {\em exponential} in $q$, however 
note that as indicated in the figure the thickness of the link is {\em linear} in $q$.
This may be thought of as a ``shared distortion'' phenomenon, an example of load balancing where 
a linear amount of local distortion creates an overall exponential effect. This phenomenon
makes the question discussed below for Milnor's invariants/Massey products interesting and probably more subtle
compared to the linking number.

We now describe the relevant $2$-complexes. 
Similar complexes were studied  in the appendix of \cite{HF} and seen
to exhibit exponential distortion, but not in ${\mathbb R}^4$, but when embedded into certain non simply connected targets.
It is an open problem whether examples exhibiting exponential intrinsic distortion exist when the
ambient $4$-manifold is ${\mathbb R}^4$. The construction below concerns embeddings into ${\mathbb R}^4$, with the focus on their thickness (rather than distortion).

A brief outline of the construction is as follows. Take two (or more) copies of the $2$-complex from section \ref{n-complexes}, then any embedding into ${\mathbb R}^4$ creates a link of $2$-spheres in ${\mathbb R}^4$. For homological reasons the fundamental group of the complement of this link modulo any finite term of the (mod 2) lower central series (see below) is isomorphic to the corresponding quotient of the free group. The complex
also provides circles that link the $2$-spheres in ${\mathbb R}^4$. Then word constructions in the free group can be implemented by attaching additional $2$-cells to the complex, and for suitable choices of the $2$-cells the complex still embeds into ${\mathbb R}^4$. Generalizing the construction in section \ref{n-complexes}, one may use $2$-cells realizing commutators to create elements in an arbitrarily high term of the lower central series. This is followed by the mapping telescope construction (\ref{telescope}). 
Since the embedding into $4$-space is not fixed, the subtlety is that the fundamental group of the $2$-link in ${\mathbb R}^4$ cannot be assumed to be actually
isomorphic to the free group, as it is controlled only modulo a term of the (mod 2) lower central series.
To illustrate this problem we pose an analogous question for classical links in ${\mathbb R}^3$ at the end of this section. We now turn to the details of the construction.

Let $\overline K_0$ denote  two copies of the complex $K_0^2$ from the statement of proposition \ref{mod 2 linking}, with two vertices (one in each copy) identified. Recall that $K^2_0$ is the
$2$-skeleton of the $6$-simplex, where one $2$-cell has been removed. Consider the boundary circle of the missing $2$-cell, and the $2$-sphere spanned in $K_0$ by the $4$ vertices which are not in the circle. Denote the two copies of these circles and $2$-spheres in $\overline K_0$ by $C',C''$ and $S', S''$ respectively. According to proposition \ref{mod 2 linking}, the complex $\overline K_0$ embeds into ${\mathbb R}^4$ and for any embedding $i$, ${\rm lk}\, (i(C'), i(S'))$ and ${\rm lk}\, (i(C''), i(S''))$ are non-zero (mod 2). Pick the two vertices that are identified in the construction of $\overline K_0$ to be in the circles $C', C''$, and denote these two loops based at the common vertex by $x, y$. This common vertex will be taken as the basepoint in the fundamental group calculations below.

Now the $2$-complex $\overline K_q$ is defined to be $\overline K_0\vee ( q$ circles labeled $w_1,\ldots, w_l)\, \cup \, (q$ $2$-cells $D_1,\ldots, D_q)$ reading
off the relations
\begin{equation} \label{2-cells}
D_1:\, w_1=[x,y], \;\; D_2:\, w_2=[x,w_1], \ldots, \; D_l:\, w_q=[x,w_{q-1}].
\end{equation}

The $2$-cells $D_i$ exhibit $w_q$ as the $q$-fold commutator
\begin{equation}\label{commutator equation}
w_q\, :=\, [x,[x,[x,\ldots, [x,y]]]].
\end{equation}

A slightly different version of the construction involves $q+1$ copies of the complex
$K_0^2$, giving rise to a commutator in distinct generators 
\begin{equation} \label{distinct labels}
w'_q\, =\, [x_1,[x_2,[x_3,\ldots, [x_q,x_{q+1}]]]].
\end{equation}

Considered in the free group $F$ on generators $x,y$, the word length of $w_q, w'_q$ is exponential in $q$. Note that it is created using $q$ $2$-cells, 
with bounded local combinatorial complexity. 
(However based on the discussion below of the corresponding problem for links in $S^3$, we expect that $\overline K_q$ embeds into
${\mathbb R}^4$ with thickness polynomial in $q$.) 

Finally, define 
$\overline K_{q,l}:=\overline K_q\cup l$-fold mapping telescope as in (\ref{telescope}), attached to the curve $w_q$. Denote by $C_{q,l}$ the left-most circle
in the mapping telescope (\ref{telescope}).
It is not difficult to see that an embedding $K_0^2\subset {\mathbb R}^4$ extends to 
an embedding of $\overline K_{q,l}$. To analyze the fundamental group of the complement,
recall the mod-$2$ version of the Stallings theorem.

Given a group ${\pi}$, define its {\em mod $2$ - lower central series} by
${\pi}_0 = {\pi}$, ${\pi}_{n+1}=[{\pi}, {\pi}_n]_2$ where $[{\pi}, {\pi}_n]_2$ is the subgroup
of $\pi$ generated by all elements of the form $fgf^{-1}g^{-1}h^2$ with $f\in{\pi}$ and $g,h\in{\pi}_n$. A theorem of Stallings \cite{Stallings} states that if a group homomorphism 
$G\longrightarrow H$ induces an isomorphism on $H_1(-;{\mathbb Z}/2)$ and an epimorphism 
on $H_2(-; {\mathbb Z}/2)$ then for each $k$ it also induces an isomorphism $G/G_k\cong H/H_k$. Moreover,  $G/G_{\omega}\longrightarrow H/H_{\omega}$
is injective, where the $\omega$-term of the mod-$2$ lower central series of $\pi$  is ${\pi}_{\omega}=\cap_{k=1}^{\infty} {\pi}_k$.

Consider any embedding $i\co \overline K_{q,l}\subset {\mathbb R}^4$.
The inclusion 
$$i(C'\vee C'')\hookrightarrow {\mathbb R}^4\smallsetminus (i(S')\cup i(S''))$$ induces a homomorphism from $F$ to 
${\pi}_1({\mathbb R}^4\smallsetminus (i(S')\cup i(S'')))$ which sends the generators $x, y$ to the based loops denoted by the same letters. By Alexander duality and the Stallings theorem, for each $k$
\begin{equation} \label{Stallings isomorphism}
F/F_k\, \longrightarrow {\pi}_1({\mathbb R}^4\smallsetminus (i(S')\cup i(S'')))/({\pi}_1)_k
\end{equation}

is an isomorphism.  Switching to the integral lower central series, observe that the curve $C_{q,l}$ represents 
a  power $w_q^N$ of the commutator (\ref{commutator equation}), where $|N|>2^l$. The integral estimates on the linking number in section \ref{n-complexes} 
do not immediately extend to the analysis of elements of higher terms of the lower central series.
It seems reasonable to expect an exponential bound, however at present this is an open problem:

{\bf Question.} Is there an exponential upper bound $c^{-l}$ on the embedding thickness of 
the complexes $\overline K_{q,l}$ in ${\mathbb R}^4$ (where $c>1$ depends on $q$)?

An analogous question can be formulated for classical links in ${\mathbb R}^3$. As discussed above, there are two versions of the question: one with the commutator (\ref{commutator equation}) in two generators $x,y$ and another version (\ref{distinct labels}) with $q+1$ distinct generators. Figure 2 shows a standard example of a link
(an iterated Bing double of the Hopf link) realizing the commutator (\ref{distinct labels}).

In this example the components $l_i$ form the unlink. More generally if they form a boundary link, then there is a map to the free group given by intersections with the spanning surfaces, and the commutators $w_q, w'_q$ are exponentially long words in the free group. 
(It is interesting to note however that the link $w\cup \{l_i\}$ may be embedded into $S^3$ with thickness {\em linear}
in the number of components, as suggested by figure 2.) To proceed with the analogy to the mapping telescope construction above, take the $N=2^l$-power of the commutator $w_q$. This is reflected in the value of the corresponding Milnor $\mu$-invariant
\cite{Milnor} in non-repeating indices.

The map (\ref{Stallings isomorphism}) is automatically an isomorphism for $2$-links in ${\mathbb R}^4$. The analogous condition holds for links in $3$-space if the $\mu$-invariants  of the link $\{ l_i\}$ are assumed to be trivial. 

{\bf Question.} Let $L$ be a $q$-component link in $S^3$ which is almost homotopically trivial (i.e. every proper sublink is homotopically trivial, in the sense of Milnor).
Let $M$ be the maximum value among Milnor's $\mu$-invariants with distinct indices 
$|{\mu}_{{i_1},\ldots, {i_q}}(L)|$.
Is there a bound thickness$(L)<c_q M^{-1}$ for some constant $c_q>0 $ independent of the link $L$? 
Is there a bound on the crossing number of $L$ in terms of $M$?

The integral formula estimates in section \ref{n-complexes} (which apply more generally to two-component links
$S^k\sqcup S^{d-k-1}\subset {\mathbb R}^d$) give an affirmative answer to the question about thickness in the case of Milnor's invariant of length $2$, equal to the linking number. 
(A bound on the ropelength in terms of the linking number for links in $S^3$ was given in \cite{CKS}.) 
The question posed above asks whether there is an analogous estimate for higher Milnor's $\mu$-invariants.

\section{Distortion of expander graphs.} \label{distortion section}

In this section we discuss distortion of embeddings into Euclidean spaces,
focusing on spaces with unbounded intrinsic distortion (\ref{intrinsic distortion}).
Specifically, we give a lower bound on distortion of expander graphs of bounded degree
with respect to embeddings into ${\mathbb R}^n$ for any fixed $n\geq 3$.
Recall that a graph $\Gamma$ is an {\em ${\alpha}$-expander}
if whenever $S$ is a subset of the set of  vertices $V$ of $\Gamma$ with $|S|\leq |V|/2$, the number of edges
connecting vertices in $S$ to vertices in $V\smallsetminus S$ is at least ${\alpha}|S|$. Interesting examples 
are families of expander graphs ${\Gamma}_k$ of bounded degree and fixed ${\alpha}>0$, with the number of vertices 
$|V({\Gamma}_k)|$ going to infinity.
Such families of examples are given by random bipartite graphs, and explicit constructions are based on groups with property (T) \cite{Margulis}, see \cite{expanders}
for a survey.

We stress that in the following statement the distortion is measured for the metric space which is the entire graph (including its edges), not just the vertex set. 
We are not aware of prior results  in the literature on the distortion of graphs considered as $1$- (rather than $0$-)complexes embedded in a Euclidean space.

\begin{theorem} \label{expander theorem} \sl
Let $\Gamma$ be an ${\alpha}$-expander of degree $\leq d$ with $|V|=N$ vertices. Then its intrinsic distortion $D_3({\Gamma})$ with 
respect to embeddings into
${\mathbb R}^3$ satisfies 
\begin{equation} \label{theorem estimate equation}
D_3({\Gamma})\, > \, C({\alpha}, d) \, N^{\frac{1}{2}}.
\end{equation}
\end{theorem}

A brief outline of the proof of theorem \ref{expander theorem}  is as follows. 
If a surface $S$ in ${\mathbb R}^3$ is chosen so that it divides the vertices of $\Gamma$ 
roughly in half then since $\Gamma$ is an $\alpha$-expander there are at least ${\alpha}N/2$ edges of $\Gamma$  intersecting $S$. Fixing a normalization, this implies the distance between these intersection points in $S$ is roughly $\sim N^{-1/2}$, and this gives a bound on the number of vertices of $\Gamma$ in a neighborhood of $S$ in terms of distortion and the degree of $\Gamma$. This argument is valid for various surfaces cutting the vertices of $\Gamma$ in half, and applying it three times shows that most vertices are in a neighborhood of the triple intersection which is zero-dimensional, giving the required estimate on distortion.

{\bf Remarks.} 
1. A direct generalization of the proof gives a bound on the distortion $D_n({\Gamma})$ for embeddings into ${\mathbb R}^n$ for any $n\geq 3$:
\begin{equation} \label{higher dimensions formula}
D_n({\Gamma})\, > \, C({\alpha}, d, n) \, N^{\frac{1}{n-1}},
\end{equation}

where the factor $C$ is of the form 
\begin{equation} \label{constant equation}
C({\alpha}, d, n)= 
c\;  n^{\frac{n}{n-1}}\; {\alpha}^{\frac{n}{n-1}}\;  d^{-1}.
\end{equation}

2. J. Matou\v{s}ek informed us that he has established an upper bound \cite{Matousek1} which almost matches our lower bound (\ref{higher dimensions formula}): a graph with $N$ vertices, maximum degree $d$, and diameter $D$,
embeds in ${\mathbb R}^n$ with distortion 
$O(D(dN)^{1/(n-2)})$. (As in theorem \ref{non-recursive thm}, the distortion is measured here for the metric space which consists of the entire graph including its edges, not just the vertex set.)

3. It is worth noting that any embedding into ${\mathbb R}^3$ of an expander graph $\Gamma$ with a large number of vertices should
contain many knot and link types.
However the theorems of \cite{GG, Pardon} on distortion of knots do not immediately carry over to the underlying graph $\Gamma$
since the distance between points in a knot $x, y\in K\subset \Gamma\subset {\mathbb R^3}$ in general is shorter in $\Gamma$ than the distance 
between $x,y$ measured in the knot $K$.

{\em Proof of theorem \ref{expander theorem}.} The proof will be given for distortion in ${\mathbb R}^3$ and at the end  we will indicate the slight modifications
needed in the general case of embeddings into ${\mathbb R}^n$. Suppose to the contrary that 
\begin{equation} \label{contrary assumption}
{\delta} \, <  \, C({\alpha}, d) \, N^{\frac{1}{2}}
\end{equation} 
where $\delta$ is the distortion of some embedding $i\! :{\Gamma}\hookrightarrow {\mathbb R}^3$ and $C({\alpha}, d)$ is the constant given in (\ref{constant equation})
for $n=3$, that is 
\begin{equation}\label{C definition}
C({\alpha},d)=
c'\; {\alpha}^{3/2}\; d^{-1}.
\end{equation}

We will identify ${\Gamma}$ with its
image under $i$. 
Note that the definition (\ref{distortion}) of ${\delta}$ is invariant
under rescaling. 
It is convenient to fix the scale as follows.
Pick any point $p$ in ${\mathbb R}^3$ and choose a radius $r_p$ such that at least $N/2$ vertices of $\Gamma$ are inside the closed ball centered at $p$ of
radius $r_p$  and at least $N/2$ vertices are in the closure of the complement of this ball. 
Consider
\begin{equation} \label{scale}
\overline r \, := \, \inf_{p\in {\mathbb R}^3} r_p.
\end{equation}

The vertices of $\Gamma$ are a discrete subset of ${\mathbb R}^3$ and $\overline r >0$ (and the infimum in (\ref{scale}) is in fact a minimum).
Let $p$ be a point such that $r_p=\overline r$.
Now rescale the embedding $i$ so that $\overline r=1$ and denote
by $S$ the sphere centered at $p$ of radius $r_p=1$.

Since $\Gamma$ is an $\alpha$-expander, there are at least ${\alpha}N/2$ edges connecting the vertices inside the ball
bounded by $S$ with the vertices which are outside of the ball. Therefore there are at least  ${\alpha}N/2$ distinct edges
intersecting the sphere $S$. Since all these intersections are contained in the $2$-sphere of radius $1$, a rough estimate shows that there are at least 
${\alpha}N/8$ distinct pairs of points $\{x,y\}$ among these intersections with 
\begin{equation} \label{distance in space}
d_{{\mathbb R}^3}(x,y)<2\, ({\alpha}N)^{-1/2}.
\end{equation}

Since the distortion equals ${\delta}$, for each such pair $\{ x,y\}$ one has 
$d_{\Gamma}(x,y)<2{\delta}({\alpha}N)^{-1/2}$. Since the edges intersecting the sphere are distinct,
any path connecting $p,q$ in $\Gamma$ must pass through a vertex of $\Gamma$. 
Moreover,
since the degree of $\Gamma$ is at most $d$, there are at least ${\alpha}N/(8d)$ vertices 
in the $2{\delta}({\alpha}N)^{-1/2}$-neighborhood of the sphere $S$. 
For simplicity of notation, denote
\begin{equation} \label{nu}
{\nu}\, := \, \frac{2 \delta}{({\alpha}\, N)^{1/2}},
\end{equation}

and let ${\mathcal N}_{{\nu}}(S)$ be the
${\nu}$-neighborhood of the sphere $S$. 
The following inequality summarizes the discussion so far:
\begin{equation} \label{summary inequality}
|V\cap {\mathcal N}_{{\nu}}(S)|\, >\, \frac{{\alpha} N}{8 d}.
\end{equation}

Picking the constant $c'$  in the expression (\ref{C definition}) to 
be $10^{-5}$, it follows from the assumption (\ref{contrary assumption}) that
\begin{equation} \label{nu estimate} 
{\nu}\; <\; \frac{{\alpha}}{10^4\, d}.
\end{equation}
Recall that $0< {\alpha}<1$ and $d>1$, so the estimate (\ref{summary inequality}) states that an
$({\alpha}/8d)$-fraction of vertices $V$ of $\Gamma$ is located in a thin $\nu$-neighborhood of the unit sphere $S$.
We will next show that under the assumption (\ref{contrary assumption}), 
in fact {\em most} of the vertices of $\Gamma$ are contained in a fairly thin neighborhood of $S$:

\begin{proposition} \label{density proposition} \sl Let ${\epsilon}>0$ (more concretely one may take ${\epsilon}=.1$ for the proof of theorem 
\ref{expander theorem} in ${\mathbb R}^3$).
\begin{equation} \label{density estimate1}
{\rm Let} \; \; \,k= \frac{10\, d}{{\alpha}\, {\epsilon}^{3/2}}.\; \; \, {\rm Then} \;\;\, |V\cap {\mathcal N}_{k{\nu}}(S)|>(1-2{\epsilon}) N.
\end{equation}
\end{proposition}

Applying this proposition to ${\epsilon}=.1$ and using the inequality (\ref{nu estimate}), one has
$k{\nu}<1/10$.  The proposition therefore asserts that a subset of the vertices $V$ of $\Gamma$ of cardinality 
$>.8 N$ is contained
in  
the $1/10$-neighborhood of the unit sphere $S$.

{\em Proof of proposition \ref{density proposition}.}
Let $B_r$ denote the ball centered at $p$ of radius $r$.
The estimate (\ref{density estimate1}) will follow from the inequalities
\begin{equation} \label{density estimate2}
|V\cap B_{1-k{\nu}}|<{\epsilon}N, \;\; |V\cap ({\mathbb R}^3\smallsetminus B_{1+k{\nu}})|<{\epsilon}N,
\end{equation}

Suppose at least one of the inequalities in
(\ref{density estimate2}) does not hold, for example $|V\cap B_{1-k{\nu}}|\geq {\epsilon}N$.

Let $r$ be any radius satisfying $1-k{\nu}<r<1$, then ${\epsilon}N<|V\cap B_r|<N/2$,
so there are at least ${\alpha}{\epsilon}N$ distinct edges of $\Gamma$ intersecting
the sphere $S_r:=\partial B_r$. As in the proof of (\ref{summary inequality})
it follows that there are at least ${\alpha}{\epsilon}N/(8d)$ vertices 
in the ${\nu}/{\epsilon}^{1/2}$-neighborhood of the sphere $S_r$.

Since $k= 10 d/({\alpha} {\epsilon}^{3/2})$, there are $10 d/({\alpha}{\epsilon})$ of such disjoint ``shells'': ${\nu}/{\epsilon}^{1/2}$-neighborhoods of spheres that fit in
$B_1\smallsetminus B_{1-k{\nu}}$, so there are a total of more than
$$\frac{10\, d}{{\alpha}{\epsilon}}\, \cdot\, \frac{{\alpha}{\epsilon}N}{8d} \; > \; N$$

vertices in this region. This contradiction concludes the proof of inequalities (\ref{density estimate2})
and of proposition \ref{density proposition}. \qed

{\em Remark}. There is a sharper version of the bound (\ref{density estimate1}): denote by $f(r)$ the
cardinality $f(r)=V\cap B_r$. The assumption that $\Gamma$ is an ${\alpha}$-expander yields a discrete
version of the differential inequality
$$f'(r)>\frac{{\alpha}f(r)}{d\, {\delta}}$$

for $r<1$ and an analogous inequality for $r>1$.
This inequality implies exponential decay for $f(k{\nu})$ as a function of $k$ which may be used to get a better constant $C({\alpha}, d, n)$
in (\ref{theorem estimate equation}), however this does not give a better estimate on the exponent of $N$.

To proceed with the proof of theorem \ref{expander theorem} for embeddings into ${\mathbb R}^3$, pick ${\epsilon}=.1$
in proposition \ref{density proposition} to conclude that a subset of the vertices $V$ of $\Gamma$ of cardinality 
$>.8 N$ is located within  
the $1/10$-neighborhood of the sphere $S$.

 Now pick a point $q\in S$ and repeat the argument: consider a sphere $S'$ centered at $q$ such that at least $N/2$ vertices of $\Gamma$ are inside $S'$ 
and at least $N/2$ vertices outside. 

Because of the normalization (\ref{scale}), the radius $r'$ of $S'$  is  $\geq 1$. However since
most of the vertices are in a small neighborhood of the unit sphere $S$, it is clear that $r'<3$. Now the
argument for $S$ given above applies verbatim to $S'$, in particular the estimate on the number of points satisfying (\ref{distance in space}) is valid
since the radii of $S$, $S'$ differ only by a bounded factor ($<3$).
This implies that a subset of the vertices of $\Gamma$ of cardinality 
$>.8N$ is located within a $1/10$-neighborhood of $S'$, so at least $.6N$ vertices are in the intersection of the
two neighborhoods ${\mathcal N} _{1/10}(S)\cap {\mathcal N}_{1/10}(S')$.

The intersection $S\cap S'$ is either (i) a circle $C$ of radius $>1/2$, or (ii) 
a smaller circle, a point or is empty. (The empty intersection $S\cap S'$  is hypothetically possible,
where the $1/10$-neighborhoods of $S, S'$ still overlap.) In the first case apply the same argument as above to a sphere $S''$ centered
at a point in $S\cap S'$ so the triple intersection is at most zero-dimensional. Moreover, at least $.4N$ vertices are in the
intersection of the $1/10$-neighborhoods of the three $2$-spheres. In this case, as well as in case (ii) above, run the argument
again for a $2$-sphere centered at the intersection point of the spheres or their neighborhoods to reach a contradiction.
This completes the proof of theorem \ref{expander theorem} for embeddings into ${\mathbb R}^3$.

The changes needed in the proof in higher dimensions $n\geq 3$ are the estimate (\ref{distance in space}) on the density of points in
the unit $(n-1)$-sphere: the exponent in general is $1/(n-1)$, and the number $n$ of $(n-1)$-spheres in ${\mathbb R}^n$ 
at the end of the proof needed to cut down the dimension of the intersection down to zero.
\qed

\medskip

{\em Acknowledgments.} We would like to thank the referee for helpful suggestions which improved the paper.

\end{document}